\documentstyle[12pt]{article}
\newcommand{\la}{\label}

\newcommand{\be}{\begin{equation}}
\newcommand{\e}{\end{equation}}
\newcommand{\ba}{\begin{eqnarray}}
\newcommand{\ea}{\end{eqnarray}}

\newcommand{\f}{\frac}

\begin{document}

\title {Existence of translating solutions to the flow by  powers of  mean curvature  on unbounded domains \thanks{
Supported by Natural Science Foundation of China (NO. 10631020 and
10871061). Part of this work was completed while H.Y. Jian was
working at Australian National University, supported by ARC DP
0879422.   }}

\author{   Huai-Yu Jian  and Hong-Jie Ju\\
  Department of  Mathematics\\
 Tsinghua University, Beijing 100084, P.R.China
 }

\date {}

\maketitle

{\bf Abstract.}  In this paper, we prove  the existence of classical
solutions of the Dirichlet problem for a class of quasi-linear
elliptic equations on  unbounded domains like a cone or a U-type
domain in $R^n (n\geq2)$. This problem comes from the study of mean
curvature flow and its generalization, the flow  by  powers of  mean
curvature.
 Our approach is a modified version of the classical Perron method, where the
solutions to the minimal surface equation are used as  sub-solutions
and a family auxiliary functions are constructed as super-solutions.


\vskip 1cm

{\large Correspondence author: Huai-Yu Jian

Email:   hjian@math.tsinghua.edu.cn

Tel: 86-10-62772864 \ \ \ \ \ \ Fax: 86-10-62785847}

\newpage

\section*{1. Introduction   }
\setcounter{section}{1} \setcounter{equation}{0}

 Given a constant $\alpha >0$ and a function $\varphi \in C^0(\partial\Omega)$.
 Consider the Dirichlet problem:  \ba\la{1.1}
\hbox{div}\left(\frac{Du}{\sqrt{1+|Du|^2}}\right)&=&-\left(\f{1}{\sqrt{1+|Du|^2}}\right)^{\alpha}\quad
\hbox{in} \ \Omega ,\\
\la{1.2}u&=&\varphi \qquad\ \ \ \ \ \ \ \ \ \ \ \ \ \ \ \ \ \ \ \
\hbox{on}\  \partial\Omega , \ea where
  $\Omega $ is an unbounded  domain in $R^n (n\geq 2)$
with $C^{2,\gamma} (0<\gamma<1)$ boundary.

The motivation to study this problem comes from   the  well-known
mean curvature flow  and its generalization, {\sl $H^k$-flow}, i.e.,
the flow of hypersurfaces by powers of mean curvature. Locally, a
{\sl $H^k$-flow } of hypersurfaces in $R^{n+1}$ can be described by
the nonlinear parabolic equation, \ba\la{1.3} \frac{\partial
V}{\partial t}=\sqrt{1+|DV|^2}\left[
\hbox{div}\left(\frac{DV}{\sqrt{1+|DV|^2}}\right) \right ]^k .\ea
When $k=1$, it is   the  well-known  mean curvature flow, which has
been studied strongly since the Huisken's work in 1984. See
[5,6,9,18,21] and the references therein.

A function $u=u(x)$ is called a {\sl translating solution} to the
$H^k$-flow   if the function $V(x,t)=u(x)+t$ solves (1.3).
Equivalently, $-u$ is a solution to equation (1.1) with $\alpha
=\frac{1}{k}.$  When $k=1,$ the translating solutions play a key
role in studying the singularity    of mean curvature flows
[5,6,8,18,21,22]. Scaling the space and time variables in a proper
way near type II-singularity points on the surfaces evolved by mean
curvature vector with a mean convex initial surface,
Huisken-Sinestrari [5,6] and White [22] proved that the limit flow
can be represented as $M_t=\{ (x, u(x)+t)\in R^{n+1}: x\in R^n, t\in
R\},$  where  $-u$ is a solution to equation (1.1) with $\alpha =1.$
Therefore, the study of type II-singularity of mean curvature flow
is reduced to studying the behavior of the solutions of equation
(1.1) with $\alpha =1$. Xu-Jia Wang [21] proved that when $\alpha
=1,$ any complete strictly convex solution of (1.1) in $R^n$ is
radially symmetric for $n=2$ and constructed a non-radially
symmetric solution  on a strip region for $n\ge 2$. Sheng and Wang
[18] used a direct argument  to study the Singularity profile in
mean curvature flow, and  the stability was studied in [1] for  the
radially symmetric solution for mean curvature flow.

 For general $k>0$,    $H^k$-flow (1.3) was studied in [15,16]. It
 was found to have important applications in minimal surfaces [2] and
isoperimetric inequalities [16].  It was proved in [19] that when
the initial surfaces are mean convex compact without boundary,  the
flow (1.3) must  blow up in finite time,  and similarly as in [5,6],
the
  type II-singularity is reduced  to the
 understanding solutions of equation (1.1) for general $\alpha
>0$.

When $\Omega$ is a bounded domain,  Marquardt [14] proved that when
 $\alpha\geq1$,  there exists   a solution in
$C^0(\bar{\Omega})\cap C^2(\Omega)$ to problem (1.1)-(1.2) if
 $\partial \Omega \in C^{2,\gamma}$, $\hbox{H}_{\partial\Omega}>0$
and $|\Omega|\leq n^n\alpha_n$.

Here and below, {\sl $\hbox{H}_{\partial\Omega}$ always denotes the
mean curvature on $\partial\Omega$ with respect to the inner normal,
and  $\alpha_n$ denotes the volume of unit ball in $R^n$.}

In [4], Gui and the authors obtained   an interior gradient
estimate, a Liouville type theorem and the asymptotic behavior at
infinity of the radially symmetric solutions to (1.1).

In this article, we prove the existence of classical solutions of
problem (1.1)-(1.2) for unbounded domains $\Omega$   like U-type or
a cone in $R^n $. To be precise, we assume that $\Omega $ satisfy
the following $(\Omega1)-(\Omega4)$.

 {\bf Assumption for $\Omega :$ }

{\sl $(\Omega1) $  there exists a sequence of bounded domains
$\{\Omega _j\}$ in $R^n$
 such that

 \ \ \ \ \ \ \ $\Omega_j \subset\Omega_{j+1}\subset
\Omega$ for any $j\geq1$ and $\Omega=\cup_{j=1}^{\infty} \Omega_j;$

$(\Omega2) $  there exists a $\gamma \in (0, 1)$ such that
 each
 $\partial\Omega_j \in C^{2, \gamma }$ and  $\hbox{H}_{\partial\Omega
 _j}>0$;

$(\Omega3)$  $\hbox{dist} (0, \Omega \setminus \Omega_j) \rightarrow
\infty $ as $j\rightarrow\infty$;

$(\Omega4)$   $\hbox{H}_{\partial\Omega}>0.$}

\vskip 0.4cm
 The main results of this paper are the following two theorems.

{\bf Theorem 1.1}  {\sl Suppose that $(\Omega1) $-$(\Omega4)$ are
satisfied and there are a constant $N $ and a positive constant $M$
such that
 $$\Omega\subset C_N(M):=\{ x=(x_1,x_2,\cdots, x_n)\in R^n| \ x_1>N,x_2^2+\cdots
+x_n^2<M^2\}$$  and   $\partial\Omega\cap \partial
C_N(M)=\emptyset$. If $\alpha >0$ and $\varphi \in
C^0(\partial\Omega)$, then there exists a solution  $u \in
C^2(\Omega) \cap C^0(\bar{\Omega})$ to problem  (1.1)-(1.2).} \vskip
0.4cm

 {\bf Theorem 1.2}  {\sl Assume that $(\Omega1)
$ -$(\Omega4)$ are satisfied  and there is a constant $\theta\in (0,
\frac{\pi}{2})$ such that
 $$\Omega\subset C(\theta):=\{ x=(x_1,x_2,\cdots, x_n)\in R^n| \ x_1>0,x_2^2+\cdots
+x_n^2<(x_1\tan \theta)^2\}$$  and $\partial\Omega\cap
\partial C(\theta)=\emptyset$. If $\alpha >0$ and  $\varphi \in
C^0(\partial\Omega)$, then there exists a solution  $u \in
C^2(\Omega) \cap C^0(\bar{\Omega})$ to problem  (1.1)-(1.2) .}\\

 The paper is organized as follows. In section 2, we prove the existence of
  Dirichlet problem (1.1)-(1.2) with $\alpha>0$ on bounded domain,
  extending the main result in [14] for the case of $\alpha \ge 1 .$
Note that
  when $0<\alpha<1$, the hypothesis  (sc)  of the corresponding
theorem in [14] can not be satisfied and the techniques in [14] can
not be applied directly.   In section 3, we construct a family of
auxiliary functions which will be used as super-solutions. In
section 4, we define the lifting function so as to construct the
class of subfunctions and  prove the properties of the subfunctions
which is necessary for the proofs of Theorems 1.1 and 1.2.  Theorems
1.1 and 1.2 will be proved in section 5 by a modified version of the
classical Perron method. The interior gradient estimate for
 (1.1) derived recently by Gui and the authors in [4] plays  an important
 role.

 \section*{2. Existence for the  solutions on bounded domains }
\setcounter{section}{2} \setcounter{equation}{0}
 In this section, we prove
the existence of the Dirichlet problem (1.1)-(1.2) with $\alpha>0$
on bounded domains, which is necessary in the proofs of theorems 1.1
and 1.2. For this purpose, we need the interior gradient estimates
 for equation (1.1), which was obtained in [4] recently by Gui and the
 authors using the idea of Xu-Jia Wang [20].
\vskip 0.4cm
 {\bf Lemma 2.1} [4]  \ {\sl Suppose $u\in C^3(B_r(0))$ is a
nonnegative solution of equation (1.1), then
$$|\nabla u(0)|\leq\exp\{C_1+C_2\frac{m^2}{r^2}\},$$
where $m=\sup_{x\in B_r(0)}u(x),$\ $C_1$ and $C_2$  are constants
depending only on $n$ and $\alpha.$ }

 {\bf Lemma 2.2}  {\sl Let   $\Omega_0\subset R^n$ be a bounded
domain with $C^{2,\gamma}$ boundary for some $\gamma \in (0,1)$ and
$|\Omega_0|< n^n\alpha_n$.  Suppose that
$\hbox{H}_{\partial\Omega_0}>0$ and $\varphi\in
C^0(\partial\Omega_0)$. Then the Dirichlet problem (1.1)-(1.2) with
$\Omega _0$ instead of $\Omega $ has a unique solution $u\in
C^0(\bar{\Omega}_0)\cap C^2(\Omega_0)$.}

{\bf Proof. }   Firstly, we suppose $\varphi \in
C^{2,\gamma}(\bar{\Omega}_0)$ and prove the Dirichlet problem
(1.1)-(1.2) has a solution $u\in C^{2,\gamma}(\bar{\Omega}_0)$. This
was proved in [14] for the case of $\alpha \ge 1$, so we assume
$\alpha \in (0, 1)$ below.

 Write (1.1)-(1.2) as
 \ba\la{eq2.1} Q u:&=& a^{ij}(Du)D_{ij}u+b(Du)=0\quad in\ \Omega_0\\
\la{eq2.2}u&=&\varphi\quad on\ \partial\Omega_0\ea where
 \ba
a^{ij}(p):&=&(1+|p|^2)\delta_{ij}-p_ip_j,\nonumber\\
b(p):&=&(1+|p|^2)^{\frac{3-\alpha}{2}}.\nonumber
 \ea
 By virtue of Theorem 13.8 in [3], it suffices to prove the $C^1$-estimate
 for the solutions $u\in C^{2,\gamma}(\bar{\Omega}_0)$ of (2.1)-(2.2).

It follows from the assumption $|\Omega_0|< n^n\alpha_n$ and Theorem
10.5 in [3] that  \ba\la{eq2.3}
\sup_{\Omega_0}|u|&\leq&\sup_{\partial\Omega_0}|u|+C
diam\Omega_0\nonumber\\
&=&\sup_{\partial\Omega_0}|\varphi|+C diam\Omega_0, \ea where
constant $C$  depends only on $n$ and $\Omega_0$.

Applying Theorem 15.1 in [3], a maximum principle for the
 gradient, we can obtain
 \begin{equation}\la{eq2.4}\sup_{\Omega_0}|Du|=\sup_{\partial\Omega_0}|Du|.\end{equation}
 Therefore, we need only to estimate $\sup_{\partial\Omega_0}|Du|$,
 which will be proved by constructing global upper and lower
 barriers for $u$ as follows.

 Let $$\Gamma:=\{x \in \bar{\Omega}_0 \ |\ d(x):=dist (x,
\partial\Omega_0) <d_1\}$$ with $0<d_1<1$ which will be determined
later. Denote $m:=\sup_{\bar{\Omega}_0}|u|$ and
$a:=\sup_{\bar{\Gamma}}|\varphi|$. We want to find a function
$\psi$, such that $w^{\pm}:=\varphi\pm\psi\circ d$ are global upper
and lower barriers for $u$ and operator $Q$ in domain $\Gamma$,
i.e.,    \ba \la{eq2.5}& &w^{\pm}=u\quad \hbox{on}\
\partial\Omega_0,\\
\la{eq2.6}& &w^-\leq u \leq w^+\quad \hbox{on}\
\partial\Gamma\setminus\partial\Omega_0,\\
\la{eq2.7}& &\pm Q w^{\pm}<0\quad \hbox{in}\
\Gamma\setminus\partial\Omega_0.\ea Assuming $\psi''(d)\leq0$ and
$\psi'(d)\geq\nu$ for some constant $\nu>0$  which will be
determined by $\Omega_0$, $\alpha$ and
$\|\varphi\|_{C^1(\bar{\Omega}_0)}$. For $x \in \Gamma$, there is a
$y \in
\partial\Omega_0$ such that $d(x)=|x-y|$.
Hence, \ba\la{eq2.8}
&\pm&a^{ij}(Dw^{\pm})D_{ij}w^{\pm}\nonumber\\
&=&\pm[(1+|Dw^{\pm}|^2)\delta_{ij}-D_iw^{\pm}D_jw^{\pm}][D_{ij}\varphi\pm\psi''D_idD_jd\pm\psi'D_{ij}d]\nonumber\\
&=&\pm(1+|Dw^{\pm}|^2)\sum_{i=1}^n D_{ii}\varphi\mp
D_iw^{\pm}D_jw^{\pm}D_{ij}\varphi\nonumber\\
& &+\psi''+\psi''[|Dw^{\pm}|^2-D_iw^{\pm}D_jw^{\pm}D_idD_jd]
\nonumber\\
&
&+\psi'(1+|Dw^{\pm}|^2)\sum_{i=1}^nD_{ii}d-\psi'D_iw^{\pm}D_jw^{\pm}D_{ij}d,
\ \  \forall x\in \Gamma, \ea where we have used the fact $|Dd|=1.$
Noting that $\psi'\geq\nu$ we have \ba\la{eq2.9}
&\pm&(1+|Dw^{\pm}|^2)\sum_{i=1}^n D_{ii}\varphi\mp
D_iw^{\pm}D_jw^{\pm}D_{ij}\varphi\nonumber\\
&\leq&2n^2(1+|D\varphi\pm \psi'
Dd|^2)\sup_{\bar{\Gamma}}|D^2\varphi|\nonumber\\
&\leq&2n^2(1+2|D\varphi|^2+2\psi'^2)\sup_{\bar{\Gamma}}|D^2\varphi|\nonumber\\
&\leq& [2n^2(\frac{1+2\sup_{\bar{\Gamma}}|D\varphi|^2}{\nu}+2
) \sup_{\bar{\Gamma}}|D^2\varphi|]\psi'^2\nonumber\\
&:=&c_1\psi'^2, \ \  \forall x\in \Gamma. \ea  By Schwarz's
inequality, \ba\la{eq2.10}
D_iw^{\pm}D_jw^{\pm}D_idD_jd\leq|Dw^{\pm}|^2.\ea Since
$D_idD_jdD_{ij}d=0$, then
\ba\la{eq2.11}&-&\psi'D_iw^{\pm}D_jw^{\pm}D_{ij}d\nonumber\\
&=&-\psi'(D_i\varphi
D_j\varphi+2\psi'D_idD_j\varphi)D_{ij}d\nonumber\\
&\leq&
[\sup_{\bar{\Gamma}}|D^2d|(\frac{n^2\sup_{\bar{\Gamma}}|D\varphi|^2}{\nu}+2n
\sup_{\bar{\Gamma}}|D\varphi|)]\psi'^2 \nonumber\\
&:=&c_2\psi'^2, \ \ \forall x\in \Gamma . \ea
 From  Lemma 14.17 in [3], \ba [D^2d(x)]=diag \left[
\frac{-k_1}{1-k_1d},\cdots,\frac{-k_{n-1}}{1-k_{n-1}d},0\right]
\nonumber\ea where $k_1,\cdots, k_{n-1}$ are the principal
curvatures of $\partial\Omega_0$ at $y$, then we have
\ba\sum^n_{i=1}D_{ii}d(x)\leq-(n-1)H_{\partial\Omega_0}(y)\nonumber\ea
if $d_1$ is small enough.  Since $\Omega_0$ is a bounded set with
$C^{2,\gamma}$ boundary  and $H_{\partial\Omega_0}>0$,
$H_0:=\min_{y\in\partial\Omega_0}H_{\partial\Omega_0}(y)=H_{\partial\Omega_0}(y_0)>0$
for some point $y_0$. Therefore,
\ba\la{eq2.12}\sum^n_{i=1}D_{ii}d(x)\leq-(n-1)H_0,\quad \forall\
x\in\Gamma.\ea Now, inserting (2.9)-(2.12) into (2.8), we obtain
\ba\la{eq2.13} \pm a^{ij}(Dw^{\pm})D_{ij}w^{\pm}\leq
\psi''+(c_1+c_2)\psi'^2-(n-1)H_0\psi'(1+|Dw^{\pm}|^2).\ea
On the
other hand, by the assumption $\alpha \in (0, 1) $ we have
\ba\la{eq2.14}|b(Dw^{\pm})|&=&(1+|Dw^{\pm}|^2)^{\frac{3-\alpha}{2}}\nonumber\\
&\leq&(1+|Dw^{\pm}|^2)[(\frac{1+2\sup_{\bar{\Gamma}}|D\varphi|^2}{\nu^2}+2)\psi'^2]^{\frac{1-\alpha}{2}}\nonumber\\
&=&(1+|Dw^{\pm}|^2)(\frac{1+2\sup_{\bar{\Gamma}}|D\varphi|^2}{\nu^2}+2)^{\frac{1-\alpha}{2}}\psi'^{1-\alpha}.
\ea Combining (2.13) and (2.14), we obtain
 \ba
\pm Qw^{\pm}&\leq&\psi''+(c_1+c_2)\psi'^2-(n-1)H_0\psi'(1+|Dw^{\pm}|^2)\nonumber\\
&+&(1+|Dw^{\pm}|^2)(\frac{1+2\sup_{\bar{\Gamma}}|D\varphi|^2}{\nu^2}+2)^{\frac{1-\alpha}{2}}\psi'^{1-\alpha}\nonumber\\
&=&\psi''+(c_1+c_2)\psi'^2-\psi'(1+|Dw^{\pm}|^2)\cdot\nonumber\\
&
&[(n-1)H_0-(\frac{1+2\sup_{\bar{\Gamma}}|D\varphi|^2}{\nu^2}+2)^{\frac{1-\alpha}{2}}\psi'^{-\alpha}].\nonumber
\ea Note that $\psi'\geq\nu$, $H_{0}>0$ and $\alpha\in (0, 1)$.
 Choose some  large number $\nu>0$ such that \ba &
&(n-1)H_0-(\frac{1+2\sup_{\bar{\Gamma}}|D\varphi|^2}{\nu^2}+2)^{\frac{1-\alpha}{2}}\psi'^{-\alpha}\nonumber\\
&\geq&(n-1)H_0-(\frac{1+2\sup_{\bar{\Gamma}}|D\varphi|^2}{\nu^2}+2)^{\frac{1-\alpha}{2}}\nu^{-\alpha}\nonumber\\
&>&0.\nonumber\ea Consequently,
   \ba\la{eq2.15}\pm
Q w^{\pm}<\psi''+(c_1+c_2)\psi'^2=:\psi''+c_3\psi'^2.\ea Thus,
(2.5)-(2.7) is reduced to  finding a function $\psi$ such that
$\psi''+c_3\psi'^2=0$, $\psi'(d)\geq\nu$, $\psi(d)\geq0$ for $d\in
(0, d_1)$, and $\psi(d_1)\geq m+a$.

Now choose the function $$\psi(d)=\frac{1}{c_3}\ln (1+kd),\ \ k>0.$$
Then
$$
\psi''+c_3\psi'^2=0,\ \  \psi (0)=0, \ \ \psi (d)>0, \ \ \forall\
d\in (0, d_1].$$
 Fix a small  $d_1\in (0,
\frac{1}{\nu c_3})$  and set
$$k=\frac{e^{c_3(a+m)}-1}{d_1}+\frac{\nu c_3}{1-\nu c_3d_1},$$
then
$$1+kd_1\geq e^{c_3(a+m)},\ \ k\geq \nu c_3(1+kd_1).$$
Thus
 $$\psi(d_1)=\frac{1}{c_3}\ln (1+kd)\geq a+m$$
 and
 $$\psi'(d)=\frac{k}{c_3(1+kd)}\geq\frac{k}{c_3(1+kd_1)}\geq\nu,\quad \hbox{for}\ 0<d\leq d_1.$$
 In this way, we have constructed barriers $w^{\pm}$ such that
 (2.5)-(2.7) are satisfied.

 Applying a maximum principle to
 (2.5)-(2.7) we see that
 $$w^-\leq u \leq w^+\quad \hbox{on}\
\partial\Gamma.$$  This, together with (2.5) again, implies
\ba\la{eq2.16}\sup_{\partial\Omega_0}|Du|\leq\sup_{\partial\Omega_0}|D\varphi|+\psi'(0)
=\sup_{\partial\Omega_0}|D\varphi|+\frac{k}{c_3}.\ea Combining
(2.3), (2.4) and (2.16), we have
\ba\la{eq2.17}\|u\|_{C^1(\bar{\Omega}_0)}=\sup_{\Omega_0}|u|+\sup_{\Omega_0}|Du|\leq
C,\ea where constant $C=C(n,\alpha, \Omega_0,\parallel
d\parallel_{C^2(\bar{\Gamma})},
\parallel \varphi\parallel_{C^2(\bar{\Omega}_0)})$. Hence,
by Theorem 13.8  in [3], the Dirichlet problem (1.1)-(1.2) has a
solution $u\in C^{2,\gamma}(\bar{\Omega}_0)$ with boundary value
$\varphi\in C^{2,\gamma}(\bar{\Omega}_0)$.

If $\varphi\in C^0(\partial\Omega_0)$, we choose a sequence of
functions $ \varphi_m \in C^{2,\gamma}(\bar{\Omega}_0)$ which is
bounded  in $C^0(\bar{\Omega}_0)$ and approximates $\varphi$ in
$C^0(
\partial\Omega_0).$   As above, the Dirichlet problem (1.1)-(1.2) has
solution $u_m \in C^{2,\gamma}(\bar{\Omega}_0)$ with boundary value
$\varphi_m$. Applying a comparison principle, $\{u_m\}$ converges
uniformly to some function $u\in C^0(\bar{\Omega}_0)$ with
$u=\varphi$ on $\partial\Omega_0$. The interior gradient estimates
(Lemma 2.1), interior H$\ddot{o}$lder estimate (Theorem 13.1 in [3])
and standard Schauder estimate imply   that there is  a subsequence
of  $\{u_m\}$ such that it converges to $u$ in
$C^{2,\gamma}(\bar{\Omega}_1)$ for any $\Omega _1\subset \subset
\Omega_0  $   by Arzel$\grave{a}$-Ascoli theorem. Thus, $u\in
C^0(\bar{\Omega}_0)\cap C^2(\Omega_0)$ solves (1.1)-(1.2). The
uniqueness
 follows directly from a comparison principle (Theorem 10.2 in [3]).  In this way, Lemma
2.2 has been proved.
$\hfill \Box$

\section*{3. A family of auxiliary functions}
\setcounter{section}{3} \setcounter{equation}{0}
 In this section, we will construct a family of auxiliary functions
 which will be used as supersolutions for problem (1.1)-(1.2).

 Recall the definition of $Qu$ in (2.1), namely,
$$Qu:=((1+|Du|^2)\delta_{ij}-D_iuD_ju)D_{ij}u+(1+|Du|^2)^{\f{3-\alpha}{2}}.$$
  We want to  construct a
family of functions $\{w_k\}$  and a family of sets $\{A_k\}$ which
covers the domains in Theorems 1.1 and 1.2, such that $Q w_k \leq 0$
in  $A_k$ for each $k\geq 1.$
  The construction method was introduced in [17] and was used again
  in [10,11] for the existence of the prescribed mean curvature
  equations in unbounded domains. Also see [13] for the existence of the constant mean curvature
  equations in unbounded convex domains.

Set
\begin{displaymath}
\Phi(\rho)=\left\{ {\ \ \rho^{-2},\quad\hbox{if}\  0<\rho<1 \atop
n-1,\quad\hbox{if}\  \rho\geq 1, \ \ \ \ \ \ }\right.
\end{displaymath}
and define a function $\xi $ by
$$\xi(t)=\int_t^{\infty}\f{d\rho}{\rho^3\Phi(\rho)}\quad \hbox{for} \,t>0. $$
Let $\eta$ be the inverse of $\xi$. It is easy to check that
\begin{displaymath}
\eta(\beta)=\left\{ {\f{1}{\sqrt{2(n-1)\beta}},\ \quad\hbox{if}\
0<\beta<\f{1}{2(n-1)}\ \ \   \atop
e^{-\beta+\f{1}{2(n-1)}},\quad\hbox{if}\
\f{1}{2(n-1)}\leq\beta<+\infty ,}\right.
\end{displaymath}
and$$\int_0^{\infty}\eta(\beta)d\beta<\infty.$$ For positive
constants $L, \mu, \tau$ with $\tau>L$ (which will be determined ),
we define\be\la{3.1} h(r)=h_{\mu,\tau}(r)=\int_r^{\tau}\eta\left(\mu
\ln\f{t}{L}\right)dt,\quad \hbox{for}\ r\in [L, \tau].\e Then $h$ is
a positive, monotonically decreasing function, satisfying
$$h(\tau)=0, \ \ h'(L)=-\infty, \ \
h(L)=\int_L^{\tau}\eta\left(\mu \ln\f{t}{L}\right)dt<\infty$$ and
\be\la{3.2} \f{h''}{(h')^3}=-\f{\mu}{r}\Phi(-h')\quad\hbox{for}\
r\in (L, \tau).\e Since $\eta(\beta) \rightarrow \infty$ as
$\beta\rightarrow0^+$, for any  $H^*>1$ there is a constant
$c(H^*,\eta)$ such that $\eta(\beta)\geq H^*$ for all
$0<\beta<c(H^*,\eta).$ Note that we may assume $c(H^*,\eta)$ is
decreasing in $H^*$.
 Letting $d=\f{c(H^*,\eta)}{\mu}$, we have \be
\la{3.3}|h'(r)|=\eta\left(\mu \ln\f{r}{L}\right)\geq H^* , \quad
\forall\ r\in (L,  Le^d) .\e Now set $\vec{x}_0=(x_1^0,0,\cdots,0),\
r(\vec{x})=|\vec{x}-\vec{x}_0|$, and \be\la{3.4}
w(\vec{x})=w_{\vec{x}_0}(\vec{x})=h(r(\vec{x})).\e Then for any $
 \vec{x} \in \{\vec{x}\in R^n|\  r(\vec{x})\in (L, Le^d) \},$  we have
$$Dw(\vec{x})=h'(r(\vec{x}))\f{\vec{x}-\vec{x}_0}{r(\vec{x})},\quad |Dw(\vec{x})|=|h'(r(\vec{x}))|\geq
H^*$$  and \ba  \label{3.5}
Qw&=&((1+|Dw|^2)\delta_{ij}-D_iwD_jw)D_{ij}w+(1+|Dw|^2)^{\f{3-\alpha}{2}}\nonumber\\
&=& h''+(n-1)(1+h'^2)\f{h'}{r}+(1+h'^2)^{\f{3-\alpha}{2}}\nonumber\\
&=& -\f{\mu}{r}h'^3\Phi(-h')+(n-1)(1+h'^2)\f{h'}{r}+(1+h'^2)^{\f{3-\alpha}{2}}\nonumber\\
&=&
|h'|^3\{\f{(n-1)\mu}{r}-\f{n-1}{rh'^2}-\f{n-1}{r}+\f{1}{|h'|^3}(1+h'^2)^{\f{3-\alpha}{2}}\},
\ea where we have used (3.2) and (3.3). \vskip 0.4cm

In order to construct the local super-solutions to equation (1.1),
we   distinguish two cases which correspond to the domains in
theorems 1.1 and 1.2 respectively.

 \noindent{\bf Case 1: } {\sl
$\Omega$ is inside the cylinder $C_N(M)$ as in Theorem 1.1.}

   Fix $0<\mu<1. $  Let $L=M$ and $ \tau=Me^d$, where
$d=\f{c(H^*, \eta)}{\mu}$ (which will be determined by $H^*$). Note
that for any fixed $\alpha>0,$\ba  \label {3.6}
\f{1}{{t}^3}(1+{t}^2)^{\f{3-\alpha}{2}}\rightarrow0\quad \hbox{as}\
t\rightarrow\infty .\ea    By (3.3) we can choose some large $H^*>1$
such that for all  $\vec{x}\in \{\vec{x}\in R^n| \  r(\vec{x})\in
(M, Me^d) \},$ \ba
\label{3.7}\f{1}{|h'|^3}(1+h'^2)^{\f{3-\alpha}{2}}
&\leq&\f{(n-1)(1-\mu)}{Me^d}\nonumber\\
&\leq&\f{(n-1)(1-\mu)}{r} .\ea   Replacing this inequality in (3.5)
we have proved

{\bf Claim 1} {\sl  For any $\mu \in (0, 1)$, there is  a $H^*>1$
such that $Q w(\vec{x})\leq0$ for all  $\vec{x} \in \{\vec{x}\in
R^n| \ r(\vec{x})\in (M, Me^d)\}$,  where $w$ is defined by (3.1)
and (3.4) with   $d=c(H^*,\eta)/\mu$, $L=M$ and $\tau=Me^d$.} \vskip
0.4cm
 For a sequence $\{a_k\},$  define $\vec{x}_k=(a_k,0,\cdots,0)$
and \be \la{3.8}A(\vec{x}_k)=\{\vec{x}=(x_1, x_2, \cdots , x_n)\in
C_N(M)|\ M< |\vec{x}-\vec{x}_k|<Me^d, x_1<a_k\}.\e By Lemma A.1 in
Appendix, we can find a small number $\varepsilon>0$ and a sequence
$\{a_k\}$ satisfying
$$a_1=N, \quad 0<a_{k+1}-a_k\leq\varepsilon M(e^d-1),\quad
k=1,2,\cdots $$ such that
$$\bigcup_{k=1}^{\infty} A(\vec{x}_k)=C_N(M)$$ and
$$\partial A(\vec{x}_{k+1})\bigcap\{\vec{x}\in C_N(M)|\
|\vec{x}-\vec{x}_{k+1}|=Me^d, x_1<a_{k+1}\} \subset A(\vec{x}_k).$$

 On each
domain $A(\vec{x}_k)$, we define a function $w_k$ as follows. Let
$h_k(r(\vec{x}))=h(|\vec{x}-\vec{x}_k|)$, where $h(r)$ is the
function defined by (3.1) with $L=M, \tau=Me^d.$
 Set
 \be\la{3.9}
 w_k(\vec{x})=h_k(r(\vec{x}))+(k-1)h(M) +\sup \{|\varphi(\vec{x})| \ |\ \vec{x}\in \partial\Omega, x_1\leq
 a_k\}.\e
  It follows from Claim 1 that each $w_k$ is well defined in
 $A(\vec{x}_k)$ and satisfies
  \be\la{3.10} Q w_k\leq 0 \ \ in \ \ A(\vec{x}_k).\e
   Furthermore, by the obvious properties of $h$, we see that
 \ba\la{3.11}
 w_k(\vec{x})&\leq&h(M)+(k-1)h(M)+\sup \{|\varphi(\vec{x})|\  |\ \vec{x}\in \partial\Omega, x_1\leq
 a_k\}\nonumber\\
 &\leq&h_{k+1}(r(\vec{x}))+k h(M)+\sup \{|\varphi(\vec{x})|\  |\ \vec{x}\in \partial\Omega, x_1\leq
 a_{k+1}\}\nonumber\\
 &=&w_{k+1}(\vec{x}), \ \ \forall\  \vec{x}\in A(\vec{x}_k)\cap A(\vec{x}_{k+1}),\ea
 where the $r$ in $h_{k+1}(r)$ is $|\vec{x}-\vec{x}_{k+1}|$.
 \vskip 0.5cm

 \noindent{\bf Case 2: } {\sl $\Omega$ is inside the cone $C(\theta)$
 as in Theorem 1.2.}

Recall that for any $L>0$, $0<\mu<1$ and   $H^*>1$ there is a
constant $c^*(H^*,\eta)$ such that (3.3) holds for
$d=\f{c(H^*,\eta)}{\mu}$, which  means that for any $0<d\leq
\f{c(H^*,\eta)}{\mu}$,   \be \la{3.12}|h'(r)|=\eta\left(\mu
ln\f{r}{L}\right)\geq H^*\quad \hbox{for}\  L<r\leq Le^d.\e

For a number $b>0$,  setting $L=b\sin\theta$, $\tau=be^d\sin\theta$
in (3.1) where $0<d\leq \f{c(H^*,\eta)}{\mu}$, we have obtained the
function $h$. Then  let $\vec{x}_0=(b,0,\cdots,0)$,
$r(\vec{x})=|\vec{x}-\vec{x}_0|$, and  $w(\vec{x})=h(r(\vec{x}))$.
It follows from (3.12) that for any   $d\in (0,  \f{c(H^*,
\eta)}{\mu})$,   $$|h'(r(\vec{x}))|\geq H^*, \ \ \forall \vec{x}\in
\{\vec{x}\in R^n| \  L< r(\vec{x})< Le^d \}.$$   Then as
(3.6)-(3.7), we have \ba \f{1}{|h'|^3}(1+h'^2)^{\f{3-\alpha}{2}}
&\leq&\f{(n-1)(1-\mu)}{Le^d}\nonumber\\
&\leq&\f{(n-1)(1-\mu)}{r}, \ \ \forall \vec{x}\in \{\vec{x}\in R^n|
\ L< r(\vec{x})< Le^d \}. \nonumber \ea   Hence, we have proved

{\bf Claim 2 } {\sl For any $b>0$, $0<\mu <1$ and $\theta \in (0,
\frac{\pi}{2}),$ there exists   $H^*>1$ such that for any $0<d\leq
c(H^*,\eta)/\mu$, $Q w \leq0$  for all  $\vec{x}\in \{\vec{x}\in
R^n\ |\  L< r(\vec{x})< Le^d \}$, where $w$ is defined by (3.1) and
(3.4)
 with $L=b\sin\theta$ and $\tau=b\sin\theta e^d$.}
 \vskip 0.4cm

 Since $\partial\Omega\cap\partial C(\theta)=\emptyset$, the vertex
of $C(\theta)$,  $0\not\in\partial\Omega$.  Hence, we can find a
small $b_1>0$ such that the ball centered at
$\vec{x}_1=(b_1,0,\cdots,0)$ with radius $b_1$ does not intersect
with $\Omega$. Choose a $d\ \in (0, \f{c(H^*,\eta)}{\mu})$ such that
$1-e^d\sin\theta>0$, and then take a   $\delta_0$ such that
\be\la{3.12}1<\delta_0<\f{1-\sin\theta}{1-e^d\sin\theta}.\e For
$k\geq 1$, let $b_k=\delta_0^{k-1}b_1$, $L_k=b_k\sin\theta$,
$\vec{x}_k=\delta_0^{k-1}\vec{x}_1$ and \be\la{3.13}
\tilde{A}(\vec{x}_k)=\{\vec{x}=(x_1, x_2, \cdots , x_n) \in
C(\theta)\ |\ L_k< |\vec{x}-\vec{x}_k|<L_ke^d,\ x_1<b_k,
  \}   . \e
By Lemma A.2 in Appendix, we have
 \be \la{3.14} \Omega \subset
 \bigcup_{k=1}^{\infty}\tilde{A}(\vec{x}_k),\e
 and the part of
$\partial\tilde{A}(\vec{x}_k)$, $S_k:=\{\vec{x}\in C(\theta)\ |\
|\vec{x}-\vec{x}_k|=L_k,\ x_1<b_k \}$,
 is completely covered by $\tilde{A}(\vec{x}_{k+1})$.

 On each domain $\tilde{A}(\vec{x}_k)$, we define a function
 $\tilde{w}_k$ as follows. Let $h_k(r)$ be the function defined by
 (3.1) with  $L=L_k=b_k\sin\theta$ and
 $\tau=L_ke^d=b_ke^d\sin\theta$. Namely,
 $$ h_k(r)=\int_{r}^{b_ke^d
 \sin\theta}\eta\left(\mu ln \f{t}{b_k\sin\theta}\right)dt, \ \  r\in [L_k , L_ke^d].$$
 Denote
 $$\tilde{B}_k=h_k(L_k)=\int_{b_k\sin\theta}^{b_ke^d
 \sin\theta}\eta\left(\mu ln \f{t}{b_k\sin\theta}\right)dt$$ and
 define
 \be\la{3.15}\tilde{w}_k(\vec{x})=h_k(r(\vec{x}))+\sum_{j=1}^{k-1}\tilde{B}_j+
 \sup\{|\varphi(\vec{x})|\ |\ \vec{x}\in \partial\Omega,\ x_1\leq
 b_k\},\e
 where $r(\vec{x})=|\vec{x}-\vec{x}_k|$. Then by Claim 2 we see that $\tilde{w}_k$ is well
 defined in $\tilde{A}(\vec{x}_k)$ and satisfies
 \be\la{3.16}Q\tilde{w}_k \leq0 \quad \hbox{in}\  \tilde{A}(\vec{x}_k).\e
Moreover,
 \ba\la{3.17}
 \tilde{w}_k(\vec{x})&\leq&h_k(L_k)+\sum_{j=1}^{k-1}\tilde{B}_j+
 \sup\{|\varphi(\vec{x})|\ |\ \vec{x}\in \partial\Omega,\ x_1\leq
 b_k\}\nonumber\\
 &\leq&h_{k+1}(r(\vec{x}))+\sum_{j=1}^{k}\tilde{B}_j+\sup\{|\varphi(\vec{x})|\ |\ \vec{x}\in \partial\Omega,\ x_1\leq
 b_{k+1}\}\nonumber\\
 &=&\tilde{w}_{k+1}(\vec{x})  \quad \hbox{in}\  \tilde{A}(\vec{x}_k)\cap \tilde{A}(\vec{x}_{k+1}),\ea
 where the $r$ in $h_{k+1}(r)$ is
 $|\vec{x}-\vec{x}_{k+1}|$.

\section*{4. The lifting and subfunction }
\setcounter{section}{4} \setcounter{equation}{0}

  In this section, we  define the lifting of
a function and the class of subfunctions which contains the
solutions of minimal surface equations. We show a few properties
which will be used   to prove the supreme function for all the
subfunctions is a solution to (1.1)-(1.2) in the next section.

  Let
$\Pi$ be the family of all bounded open sets $ O\subset \Omega$
satisfying  $\partial O\in C^{2,\gamma} $, $H_{\partial O}
>0$ and $|O|< n^n\alpha_n$.  $\varphi,  $ $C_N(M)$ and $ C(\theta)$
are the same as in Theorems 1.1 and 1.2.

{\bf Definition 4.1}  {\sl Let $v \in C^0(\bar{\Omega})$.  For each
$O \in \Pi$, define a new function $M_O(v)$, called the lifting of
$v$ over $O$, as follows:
\begin{displaymath}
M_O(v)(\vec{x})=\left\{ {v(\vec{x}),\quad\hbox{if}\  \vec{x}\ \in
\Omega\setminus O \atop z(\vec{x}),\quad\hbox{if}\  \vec{x}\ \in O\
\ \ \ \ \ }\right.
\end{displaymath}
where $z(\vec{x})$ is the solution of the boundary-value problem
\begin{displaymath}
\left\{ {Q z=0,\quad\hbox{in}\  O,\ \ \ \ \ \ \ \   \atop\ \
z=v,\quad\hbox{on}\
\partial O.\ \ \ \ \ }\right.
\end{displaymath}}
Note that the definition is well-defined by Lemma 2.2.

{\bf Definition 4.2} {\sl  The subfunction class $\Large F$ is
defined as follows: a function $v$ is in $\Large F$ if and only if

(1) $v \in C^0(\bar{\Omega})$ and $v\leq\varphi$ on
$\partial\Omega;$

(2) for any $O \in \Pi$, $v\leq M_O(v)$;

(3) if $\Omega\subseteq C_N(M)$, then $v\leq w_k$ in $\Omega \cap
A(\vec{x}_k)$ for $k\geq 1$;

(4) if $\Omega\subseteq C(\theta)$, then $v\leq \tilde{w}_k$ in
$\Omega \cap \tilde{A}(\vec{x}_k)$ for $k\geq 1$.} \vskip 0.4cm

 Let
$\Omega_1$ be a domain in $R^n$, $\{c_k\}_{k=1}^{\infty}$ be a
non-negative, non-decreasing sequence. If $\Omega_1$ is inside the
cylinder $C_N(M)$, we set \be\la{4.1}
w_k^1(\vec{x})=h_k(r(\vec{x}))+(k-1)h(M) +c_k\quad \hbox{in}\
A(\vec{x}_k),\e where $h_k$ and $A(\vec{x}_k)$ are the same as those
defined in (3.8) and (3.9). Thus, $w_k^1$ satisfies (3.10) and
(3.11) in $A(\vec{x}_k)$.

 If $\Omega_1$ is inside the cone $C(\theta)$, we set
 \be\la{4.2}
\tilde{w}_k^1(\vec{x})=h_k(r(\vec{x}))+\sum_{j=1}^{k-1}\tilde{B}_j
+c_k\quad \hbox{in}\ \tilde{A}(\vec{x}_k),\e where $h_k$,
$\tilde{B}_j$ and $\tilde{A}(\vec{x}_k)$ are the same as   in (3.14)
and (3.16).  Thus, $\tilde{w}_k^1$ satisfies (3.17) and (3.18) in
$\tilde{A}(\vec{x}_k)$.

{\bf Lemma 4.1} {\sl  Suppose  $u \in C^2(\Omega_1)\cap
C^0(\bar{\Omega}_1)$  and $Qu\geq 0$ in  $\Omega_1.$

 (i) \ When $\Omega_1\subset C_N(M) $ and  $\partial\Omega_1 \cap \partial
C_N(M)=\emptyset$,   if \be\la{4.3} u\leq w_k^1 \quad \hbox{on}\
A(\vec{x}_k)\cap\partial\Omega_1\quad \hbox{for} \ k\geq1,\e
 then $$u\leq
w_k^1 \quad \hbox{in}\ A(\vec{x}_k)\cap\Omega_1\quad \hbox{for} \
k\geq1.$$

(ii) \ When $\Omega_1\subset C(\theta) $ and $\partial\Omega_1 \cap
\partial C(\theta)=\emptyset$, if \be\la{4.4} u\leq \tilde{w}_k^1 \quad \hbox{on}\
\tilde{A}(\vec{x}_k)\cap\partial\Omega_1\quad \hbox{for} \ k\geq1,\e
 then $$u\leq
\tilde{w}_k^1 \quad \hbox{in}\ \tilde{A}(\vec{x}_k)\cap\Omega_1\quad
\hbox{for} \ k\geq1.$$ }

{\bf Proof} \  At first, let us  prove (i).

Among the family of domains $A(\vec{x}_k)$, let $A(\vec{x}_{k_0})$
be the first one (i.e. smallest $k$) which  intersects with
$\Omega_1$. We  conclude that \be\la{4.5} u\leq w_{k_0}^1\quad
\hbox{in}\ A(\vec{x}_{k_0})\cap\Omega_1.\e In fact, by (4.3),
$$u\leq w_{k_0}^1 \qquad \hbox{on}\  A(\vec{x}_{k_0})\cap\partial\Omega_1.$$
Note that $\partial A(\vec{x}_{k_0})\cap\Omega_1\cap\{
|\vec{x}-\vec{x}_{k_0}|=Me^d\}$ is empty.  Otherwise, from the fact
that $\partial A(\vec{x}_{k_0})\cap\{|\vec{x}-\vec{x}_{k_0}|=Me^d\}$
is covered by $A(\vec{x}_{k_0-1})$ (Lemma A.1),  we see that
$A(\vec{x}_{k_0})$ will not be the first to intersect with
$\Omega_1$,  a contradiction.

Also, $\partial
A(\vec{x}_{k_0})\cap\Omega_1\cap\{M<|\vec{x}-\vec{x}_{k_0}|<Me^d\}$
is empty, which follows from the fact that $\partial
A(\vec{x}_{k_0})\cap\{M<|\vec{x}-\vec{x}_{k_0}|<Me^d\}$ is a part of
$\partial C_N(M)$ and $\partial\Omega_1\cap\partial
C_N(M)=\emptyset$ by the assumption.

On $\partial
A(\vec{x}_{k_0})\cap\Omega_1\cap\{|\vec{x}-\vec{x}_{k_0}|=M\}$, it
follows from the fact $h'(M)=-\infty $ that the outer normal
derivative of $w_{k_0}^1$ is $+\infty$. Thus, $u-w_{k_0}^1$ cannot
achieve a maximum on this part of the boundary.

Therefore, $$u\leq w_{k_0}^1\ \ \hbox{on}\ \partial(
A(\vec{x}_{k_0})\cap\Omega_1).$$ Furthermore, (3.10) and the
assumption imply
$$Q w_{k_0}^1\leq Q u
\quad \hbox{in}\ A(\vec{x}_{k_0})\cap\Omega_1.$$  Hence (4.5)
follows from   the standard maximum principle [3].

 We now compare
$u$ with $w_{k_0+1}^1$ on $A(\vec{x}_{k_0+1})\cap\Omega_1.$  By
(4.3),
$$u\leq w_{k_0+1}^1 \quad \hbox{on}\
A(\vec{x}_{k_0+1})\cap\partial\Omega_1.$$
 Since $\partial
A(\vec{x}_{k_0+1})\cap\Omega_1\cap\{|\vec{x}-\vec{x}_{k_0+1}|=Me^d\}$
is covered by $A(\vec{x}_{k_0})$(Lemma A.1), then $u\leq
w_{k_0}^1\leq w_{k_0+1}^1$ on this part, by   (4.5) and (3.11).

As above, $\partial
A(\vec{x}_{k_0+1})\cap\Omega_1\cap\{M<|\vec{x}-\vec{x}_{k_0+1}|<Me^d\}$
is also empty.

On $\partial
A(\vec{x}_{k_0+1})\cap\Omega_1\cap\{|\vec{x}-\vec{x}_{k_0+1}|=M\}$,
the outer normal derivative of $w_{k_0+1}^1$ is $+\infty$. Thus,
$u-w_{k_0+1}^1$ cannot achieve a maximum on this part of the
boundary.

Since $$Qw_{k_0+1}^1\leq Q u \quad \hbox{in}\
A(\vec{x}_{k_0+1})\cap\Omega_1,$$  by the standard maximum principle
[3] we obtain \be\la{4.6}u-w_{k_0+1}^1\leq0\quad\hbox{in}\
A(\vec{x}_{k_0+1})\cap\Omega_1.\e
Repeating the above procedure, we
can obtain $$u\leq w_k^1 \quad \hbox{in}\ A(\vec{x}_k)\cap\Omega_1,\
\ \forall k\geq 1.
$$

The proof of (ii) is almost the same, and we write as follows just
for the completeness.  In the family of domains
$\tilde{A}(\vec{x}_k)$, let $\tilde{A}(\vec{x}_{k_0})$ be the first
one (i.e.smallest $k$) to intersect with $\Omega_1$. We first
conclude that \be\la{4.7} u\leq \tilde{w}_{k_0}^1 \quad \hbox{in}\
\tilde{A}(\vec{x}_{k_0})\cap\Omega_1.\e In fact, by (4.4) we have
$$u\leq \tilde{w}_{k_0}^1 \qquad \hbox{on}\  \tilde{A}(\vec{x}_{k_0})\cap\partial\Omega_1.$$

Note that $\partial
\tilde{A}(\vec{x}_{k_0})\cap\Omega_1\cap\{|\vec{x}-\vec{x}_{k_0}|=L_{k_0}e^d\}$
is empty.  Otherwise, by  the fact that $\partial
\tilde{A}(\vec{x}_{k_0})\cap\{|\vec{x}-\vec{x}_{k_0}|=L_{k_0}e^d\}$
is covered by $\tilde{A}(\vec{x}_{k_0-1})$ (Lemma A.2),  we see that
$\tilde{A}(\vec{x}_{k_0})$ will not be the first to intersect with
$\Omega_1$, a contradiction.

$\partial
\tilde{A}(\vec{x}_{k_0})\cap\Omega_1\cap\{L_{k_0}<|\vec{x}-\vec{x}_{k_0}|<L_{k_0}e^d\}$
is also empty, since $\partial
\tilde{A}(\vec{x}_{k_0})\cap\{L_{k_0}<|\vec{x}-\vec{x}_{k_0}|<L_{k_0}e^d\}$
is a part of $\partial C(\theta)$ and $\partial\Omega_1\cap\partial
C(\theta)=\emptyset$ by  the assumption.

On $\partial
\tilde{A}(\vec{x}_{k_0})\cap\Omega_1\cap\{|\vec{x}-\vec{x}_{k_0}|=L_{k_0}\}$,
the outer normal derivative of $\tilde{w}_{k_0}^1$ is $+\infty$.
Thus, $u-\tilde{w}_{k_0}^1$ cannot achieve a maximum on this part of
the boundary.

Since $$Q\tilde{w}_{k_0}^1\leq Qu \quad \hbox{in}\
\tilde{A}(\vec{x}_{k_0})\cap\Omega_1,$$ by a maximum principle we
obtain \be\la{4.8}u-\tilde{w}_{k_0}^1\leq0\quad\hbox{in}\
\tilde{A}(\vec{x}_{k_0})\cap\Omega_1.\e

We now compare $u$ with $\tilde{w}_{k_0+1}^1$ in
$\tilde{A}(\vec{x}_{k_0+1})\cap\Omega.$
 By (4.4) again, $$u\leq
\tilde{w}_{k_0+1}^1 \quad \hbox{on}\
\tilde{A}(\vec{x}_{k_0+1})\cap\partial\Omega_1.$$
 Since $\partial
\tilde{A}(\vec{x}_{k_0+1})\cap\Omega_1\cap\{|\vec{x}-\vec{x}_{k_0+1}|=L_{k_0+1}e^d\}$
is covered by $\tilde{A}(\vec{x}_{k_0})$ (Lemma A.2), then $u\leq
\tilde{w}_{k_0}^1\leq \tilde{w}_{k_0+1}^1$ on this part, by (4.8)
and (3.18).

As above, $\partial
\tilde{A}(\vec{x}_{k_0+1})\cap\Omega_1\cap\{L_{k_0+1}<|\vec{x}-\vec{x}_{k_0+1}|<L_{k_0+1}e^d\}$
is also empty.

On $\partial
\tilde{A}(\vec{x}_{k_0+1})\cap\Omega_1\cap\{|\vec{x}-\vec{x}_{k_0+1}|=L_{k_0+1}\}$,
the outer normal derivative of $\tilde{w}_{k_0+1}^1$ is $+\infty$.
Thus, $u-\tilde{w}_{k_0+1}^1$ cannot achieve a maximum on this part
of the boundary.

Since $$Q\tilde{w}_{k_0+1}^1\leq Qu \quad \hbox{in}\
\tilde{A}(\vec{x}_{k_0+1})\cap\Omega_1,$$ by a maximum principle we
obtain \be\la{4.9}u-\tilde{w}_{k_0+1}^1\leq0\quad\hbox{in}\
\tilde{A}(\vec{x}_{k_0+1})\cap\Omega_1.\e  Repeating the above
procedure as necessary,   we   arrive at $$u\leq \tilde{w}_k^1 \quad
\hbox{in}\ \tilde{A}(\vec{x}_k)\cap\Omega_1, \ \ \forall k\geq 1. $$
$\hfill \Box$

 {\bf Corollary
4.1} {\sl Let $\Omega $ be the same domain as in Theorem 1.1. If $u
\in C^2(\Omega)\cap C^0(\bar{\Omega})$ be a solution of the  problem
\ba\la{4.10}((1+|Du|^2)\delta_{ij}-D_iuD_ju)D_{ij}u&=&0 \ \quad
\hbox{in} \
\Omega\\
 \la{4.11}u&=&\varphi\quad \hbox{on}\  \partial\Omega,\ea
 then $$|u(\vec{x})|\leq w_k(\vec{x}) \quad \hbox{in}\
A(\vec{x}_k)\cap\Omega,\  \hbox{for}\  k\geq1.$$
 }

{\bf Proof}      \  Note that
$$Q u\geq0\quad \hbox{in} \  A(\vec{x}_k)\cap\Omega,$$
$$Qw_k\leq0\quad \hbox{in} \  A(\vec{x}_k)\cap\Omega$$
and  $u=\varphi\leq\sup\{|\varphi(\vec{x})|\ |\ \vec{x} \in
\partial\Omega,\ x_1<a_k\}\leq w_k(\vec{x})$ on
$A(\vec{x}_k)\cap\partial\Omega$. By the conclusion (i) of  Lemma
4.1, we can obtain $$u\leq w_k\quad \hbox{in} \
A(\vec{x}_k)\cap\Omega,\ k\geq1.$$
On the other hand, $v=-u$ is also
a solution of (4.10)-(4.11) with $\varphi$  replaced by $-\varphi$,
we can get$$v=-u\leq w_k\quad \hbox{in} \  A(\vec{x}_k)\cap\Omega,\
k\geq1.$$ Therefore,
$$|u|\leq w_k\quad \hbox{in} \  A(\vec{x}_k)\cap\Omega,\
k\geq1.$$

$\hfill \Box$

 Similarly, we have

{\bf Corollary 4.2} {\sl Let $\Omega $ be the same domain as in
Theorem 1.2.  If $u \in C^2(\Omega)\cap C^0(\bar{\Omega})$ be a
solution of the problem
\ba((1+|Du|^2)\delta_{ij}-D_iuD_ju)D_{ij}u&=&0 \quad\ \hbox{in} \
\Omega\nonumber\\
 u&=&\varphi\quad \hbox{on}\  \partial\Omega,\nonumber\ea
 then $$|u(\vec{x})|\leq \tilde{w}_k\quad \hbox{in}\
\tilde{A}(\vec{x}_k)\cap\Omega,\  \hbox{for}\  k\geq1.$$
 }

{\bf Corollary 4.3} {\sl Let $\Omega $ be the same domains as in
Theorems 1.1 or  1.2. Then $\Large F$ is not empty. }

{\bf Proof}  It follows from Lemma 4.5 in [9] that  under the
assumption $(\Omega1)-(\Omega3)$,  the boundary-value problem
(4.10)-(4.11) has a solution $v_0 \in C^2(\Omega)\cap
C^0(\bar{\Omega})$.  By Corollaries 4.1 or 4.2, we can see that $v_0
\in \Large F.$

$\hfill \Box$

 \vskip 0.4cm
 Next, we show a few properties of subfunctions which will be
 necessary in the proofs of theorems 1.1 and 1.2.
 For this purpose,
we assume that $\Omega $  is one of the following cases:

{\bf Case (i)}\  {\sl $\Omega\subset C_N(M) $ and $\partial\Omega
\cap \partial C_N(M)=\emptyset$;  }

{\bf Case (ii)}\  {\sl $\Omega\subset C(\theta) $ and
$\partial\Omega\cap \partial C(\theta)=\emptyset$.}

First, we assume case (i) and prove the following three lemmas,
which also hold for case (ii).

\vskip 0.4cm {\bf Lemma 4.2}  {\sl If $v_1, v_2 \in
C^{0}(\bar{\Omega})$ and $v_1\leq v_2$ in $\Omega $, then
$M_O(v_1)\leq M_O(v_2)$ for any $O \in \Pi.$}

 {\bf Proof} \ \  By the definition of $M_O(v)$, we have $$M_O(v_1)=v_1\leq v_2=M_O(v_2)\quad \hbox{on}\ \Omega\setminus
 O,$$ thus, we need only to prove $M_O(v_1)\leq M_O(v_2)$ on $O$.\\
 Since $z_i:=M_O(v_i)(i=1,2)$ satisfies the boundary-value equation
\ba Q z_i &=&0\quad
\hbox{in}\ O\nonumber\\
z_i&=&v_i\ \ \hbox{on}\ \partial O\nonumber\ea   and $z_1=v_1\leq
v_2=z_2$ on $\partial O$, then by a comparison principle we obtain
$z_1\leq z_2$ in $O$. Therefore, $M_O(v_1)\leq M_O(v_2)$ on
$\Omega$. $\hfill \Box$

\vskip 0.4cm {\bf Lemma 4.3}  {\sl If $v_i \in F \ \  (i=1, 2)$,
then $\max\{v_1,v_2\} \in \Large F.$}

{\bf Proof}\  \ \  By the definition of $\Large F$,   $\max\{v_1,
v_2\}\in C^0(\bar{\Omega})$, $\max\{v_1, v_2\}\leq\varphi$ on
$\partial\Omega$ and $\max\{v_1, v_2\}\leq w_k$ in $A(x_k)\cap
\Omega$ for $k\geq1$. So we  need only to check that
for any $O \in \Pi$, $\max\{v_1,v_2\}\leq M_O(\max\{v_1,v_2\})$.\\
Since $v_i\leq \max\{v_1,v_2\} \ \  (i=1, 2)$, by Lemma 4.2 we have
that for any $O \in \Pi$, $$M_O(v_i)\leq M_O(\max\{v_1,v_2\}) \ \ \
(i=1, 2) .$$ Since $v_i \in F$ imply that $v_i\leq M_O(v_i)$,
  we obtain $v_i\leq M_O(\max\{v_1,v_2\}) \ \  (i=1, 2).$
  Namely, $\max\{v_1,v_2\}\leq M_O(\max\{v_1,v_2\})$.
 $\hfill
\Box$

\vskip 0.4cm  {\bf Lemma 4.4}  {\sl If $v \in \Large F$, then
$M_O(v) \in \Large F$ for any $O \in \Pi.$}

{\bf Proof}\ \  By the definition, $M_O(v)\in C^0(\bar{\Omega})$ and
$M_O(v)=v\leq \varphi$ on $\partial\Omega$ .

First we prove that, for any $O_1\in\Pi$, \be \la{4.12}M_O(v)\leq
M_{O_1}(M_O(v)).\e Observe that
\be\la{4.13}M_{O_1}(M_O(v))=M_O(v)\quad \hbox{in}\ \Omega\setminus
O_1.\e  It is enough to prove that (4.12) holds on $O_1$.

 Since $v\leq M_O(v)$ on $\Omega$, then we have
$M_{O_1}(v)\leq M_{O_1}(M_O(v)) $ by Lemma 4.2. Moreover, we have
\be\la{4.14}M_{O}(v)=v\leq M_{O_1}(M_O(v)) \quad \hbox{in}\
O_1\setminus O .\e Denote $z_1=M_O(v)$ and $z_2=M_{O_1}(M_O(v))$. We
see that
$$ Qz_i =0\quad \hbox{in}\ O_1\cap O, \ \ i=1,2.$$
It follows from (4.13), (4.14) and the continuity of $z_i $ that
\be\la{4.15} z_1=M_{O}(v)\leq M_{O_1}(M_O(v))=z_2\quad \hbox{on}\
\partial(O_1\cap O).\e
Then a comparison  principle implies that $z_1\leq z_2$ in $O_1\cap
O$. Thus, (4.12) is true in $O_1\cap O$ and hence in $O_1$ by
(4.14).

It remains  to prove that $M_O(v)\leq w_k$ in $A(\vec{x}_k)\cap
\Omega$ for all $k\geq1$.   Since $v \in \Large F$, we find that \ba
M_O(v)=v\leq w_k\quad \hbox{in}\ A(\vec{x}_k)\cap\partial O , \ \
\forall k\geq1.\nonumber\ea Thus, the assumption (4.3) in Lemma 4.1
is satisfied for  $\Omega_1=O$, $w_k^1=w_k $ and $c_k=\sup
\{|\varphi(\vec{x})| \ |\ \vec{x}\in
\partial\Omega, x_1\leq
 a_k\} $.  Apply  this lemma to $u=M_O(v)$ we  conclude
 that $M_O(v)\leq w_k$  in $A(\vec{x}_k)\cap O.$

$\hfill \Box$

 If  case  (ii) happens,  replacing $w_k$ and $A(\vec{x}_k)$ by $\tilde{w}_k$ and
$\tilde{A}(\vec{x}_k)$, respectively, without changing the rest of
the proof, we  see that Lemmas 4.2, 4.3 and 4.4 also hold.

\section*{5.  Proofs of Theorems 1.1 and 1.2}
\setcounter{section}{5} \setcounter{equation}{0}

 We are in the position to use Perron's method to prove the
 theorems.

{\bf Proof of Theorem 1.1:}  Set $u(\vec{x})=\sup\{v(\vec{x})\ |\ v
\in \Large F\}$ for  $ \vec{x}\in \bar{\Omega}.$  We will show that
$u$ is in $C^0(\bar{\Omega})\cap C^2(\Omega)$ and satisfies
(1.1)-(1.2).

For any $\vec{x}_0 \in \Omega$,  by the definition of
$u(\vec{x}_0)$, there is a sequence of functions
$\{v_i\}_{i=1}^{\infty} \subset  \Large F$ such that
$$u(\vec{x}_0)=\lim_{i\rightarrow\infty}v_i(\vec{x}_0).$$
Let $v_0$ be a solution of (4.10)-(4.11). Then by the proof of
Corollary 4.3, we have
 \be \la {5.1} v_0 \in {\Large F }\ \ \ and \ \ \ u\ge v_0 \ \ in \ \
 \Omega. \e
 Replacing $v_i$ by $\max\{v_i, \ v_0\}$, we may assume that $v_i\geq
v_0$ on $\Omega$ by Lemma 4.3. For any $O \in \Pi$ such that
$\vec{x}_0 \in O$, replacing $v_i$ by $M_O(v_i)$, we then obtain a
sequence of functions $z_i=M_O(v_i)$  such that
$$u(\vec{x}_0)=\lim_{i\rightarrow\infty}z_i(\vec{x}_0),$$
\ba
Q z_i &=&0\quad \hbox{in} \ O , \nonumber\\
z_i&=&v_i\ \ \hbox{on}\ \partial O . \nonumber\ea Since, for all $k
$ and  $i$, \be \la {5.2} v_0\leq v_i\leq z_i\leq w_k\quad
\hbox{in}\ O\cap A(\vec{x}_k)\e  and  $O$ can be covered by the
finitely many   domains $A(\vec{x}_k)$, there is a  constant $K_1$
  such that
$$v_0\leq z_i\leq K_1\quad \hbox{in}\ O, \ \  \forall i\geq 1 .$$
Using Lemma 2.1 first, then  the standard interior
H$\ddot{\hbox{o}}$lder estimate of the gradients [3, Theorem 13.1]
and finally standard Schauder estimates [3], by
Arzel$\grave{\hbox{a}}$-Ascoli theorem we can choose    a
subsequence of $z_i$ ( denoted still by $z_i$) converging to a
function $z \in C^2(O)$ and so  $z(\vec{x})$ satisfies \be \la
{5.3}Q z=0\quad \hbox{in} \ O.\nonumber\e Obviously,
$u(\vec{x}_0)=z(\vec{x}_0)$ and $u(\vec{x})\geq z(\vec{x})$ in $O$.

Next, we prove that $u\equiv z$ on $O$. Indeed, if there is another
point $\vec{x}_1 \in O$ such that $u(\vec{x}_1)>z(\vec{x}_1)$, then
there is a function $u_0\in \Large F$ such that
$$z(\vec{x}_1)<u_0(\vec{x}_1)\leq u(\vec{x}_1).$$ Setting
$\bar{z_i}=M_O(\max\{u_0,M_O(v_i)\}),$   we have,  for all $k $ and
$i$, that
$$v_0\leq v_i\leq \bar{z}_i\leq w_k\quad \hbox{in}\ O\cap
A(\vec{x}_k)$$ and $Q \bar{z_i}=0$ in $O$.  Repeating the arguments
from (5.2) to (5.3), we obtain a subsequence of $\{\bar{z}_i\}$ (
denoted still by $\bar{z}_i$) which converges to a
  function $\bar{z}$ in $C^2(O)$  and $Q\bar{z}=0$ on $O$.
Obviously $$z_i=M_O(v_i)\leq M_O(\max\{u_0,M_O(v_i)\})=\bar{z}_i.
$$
 Hence, $$z\leq \bar{z}\quad \hbox{in}\ O,$$
$$z(\vec{x}_1)<u_0(\vec{x}_1)\leq \bar{z}(\vec{x}_1)$$and$$z(\vec{x}_0)=u_0(\vec{x}_0)=\bar{z}(\vec{x}_0).$$
That is, $\bar{z}(\vec{x})-z(\vec{x})$ is non-negative, not
identically zero in $O$ and attains its minimum value zero inside
$O$. However, it follows from the equations satisfied by $z$ and
$\bar{z}$, we find that \ba
& &((1+|D\bar{z}|^2)\delta_{pq}-D_p\bar{z}D_q\bar{z})D_{pq}(\bar{z}-z)\nonumber\\
&=&E(x,z,\bar{z},Dz,D\bar{z},D^2z,D^2\bar{z})D(\bar{z}-z)\quad
\hbox{in} \ O\nonumber\ea for some continuous function $E$. Then, by
the standard maximum principle, we have got a contradiction. Thus,
$u\equiv z$ in $O$. Since $O$ can be arbitrary, $u\in C^2(\Omega)$
and $Q u =0$ in $\Omega$.

Finally, it remains to  prove that $$u\in C^0(\bar{\Omega})\ \ and \
\ u=\varphi \ \ on \ \ \partial \Omega.$$
 For any point $\vec{x}_2
\in
\partial\Omega$, we can find a bounded $C^{2,\gamma}$ domain
$\Omega_1\subset\Omega$ such that
 $\partial\Omega_1 \cap
\partial \Omega$ is an open neighborhood of $\vec{x}_2$ in
$\partial\Omega$, $|\Omega_1|< n^n\alpha_n$  and
 $H_{\partial\Omega_1}>0.$

Since $\Omega_1$ is covered by finitely many $A(\vec{x}_k)$, there
is a  constant $K_3>0$ such that
 \be \la {5.4} v\leq K_3\quad \hbox{in}\
\bar{\Omega}_1, \ \  \forall v \in F .\e
 Now on $\partial\Omega_1$, we
choose a continuous function $\varphi^*$ as follows: $\varphi^*=K_3$
on $\partial\Omega_1 \cap\Omega$; $\varphi^*=\varphi$ in a
neighbourhood of $\vec{x}_2$ in $\partial\Omega_1 \cap \partial
\Omega$; and $\varphi^*\geq\varphi$ on the rest of $\partial\Omega_1
$. Consider the boundary value problem \ba \la{5.5}Q u &=&0\quad\
\hbox{in} \ \Omega_1,\\
\la{5.6}u&=&\varphi^*\ \  \hbox{on}\ \partial\Omega_1,\ea  which has
a solution $u_1\in C^2(\Omega_1)\cap C^0(\bar{\Omega}_1)$  by Lemma
2.2.  Therefore, for any $v \in \Large F$ we have $M_O(v)\leq u_1$
in  $ \Omega_1.$ Hence, $u\leq u_1$ in  $\Omega_1,$ which together
with (5.1), implies
$$v_0\leq u\leq u_1\quad \hbox{on}\ \Omega_1.$$ The continuity of
$u$ at $\vec{x}_2$ then follows from the fact that $v_0=u_1=\varphi$
on a neighbourhood of $\vec{x}_2$ in $\partial\Omega$ and both $v_0$
and $u_1$ are continuous in a neighbourhood of $\vec{x}_2$ in
$\bar{\Omega}$. Since $\vec{x}_2 \in
\partial\Omega$ can be arbitrary, we have proved $u \in
C^0(\bar{\Omega})$ and $ u=\varphi $ on $ \partial \Omega .$

$\hfill \Box$

{\bf Proof of Theorem 1.2: } In   this case,  $\Omega$ is inside $
C(\theta)$.  Replacing $w_k$ and $A(\vec{x}_k)$ by $\tilde{w}_k$ and
$\tilde{A}(\vec{x}_k)$ respectively, without changing the rest of
the proof of Theorem 1.1, we can obtain Theorem 1.2. $\hfill \Box$

\section*{Appendix A}
\setcounter{section}{6} \setcounter{equation}{0}

 {\bf Lemma A.1}
{\sl Let $M, d$ be positive constant. For a sequence $\{a_k\}$, set
$\vec{x}_k=(a_k,0,\cdots,0)$ and
$$A(\vec{x}_k)=\{\vec{x}=(x_1, \cdots ,  x_n)\in C_N(M)|\ M< |\vec{x}-\vec{x}_k|<Me^d,
x_1<a_k\}  . $$ Then there exists a $\varepsilon \in (0,1)$ such
that if $\{a_k\}$ satisfies
 \be\la{6.1}a_1=N, \quad 0<a_{k+1}-a_k\leq\varepsilon
M(e^d-1),\quad k=1,2,\cdots£¬\e
  then the part of
the boundary of $A(\vec{x}_{k+1})$
$$\{\vec{x}=(x_1, \cdots ,  x_n)\in C_N(M)|\  |\vec{x}-\vec{x}_{k+1}|=Me^d, x_1<a_{k+1}\}$$is inside
$A(\vec{x}_k)$. Thus, $$C_N(M)=\bigcup_{k}A(\vec{x}_k).$$ }

 {\bf
Proof}\ \
 For $\vec{x}\in \{\vec{x}\in C_N(M)|\
 |\vec{x}-\vec{x}_{k+1}|=Me^d, x_1<a_{k+1}\}$, we have \be\la{6.2}
(x_1-a_{k+1})^2+\sum_{i=2}^nx_i^2=M^2e^{2d} \e and
\be\la{6.3}x_1<a_{k+1}.\e We need only to prove that
\be\la{6.4}M^2<(x_1-a_{k})^2+\sum_{i=2}^nx_i^2<M^2e^{2d}\e and
\be\la{6.5}x_1\leq a_k.\e We first   verify (6.5).  In fact, by
(6.2) and the definition of $C_N(M)$, we have \ba (x_1-a_{k+1})^2
>M^2e^{2d}-M^2,\nonumber\ea which, together with (6.3),
implies
 \ba x_1<a_{k+1}-M\sqrt{e^{2d}-1}.\nonumber\ea In order to
prove (6.5), it is sufficient to show that
 \be\la{6.6}
a_{k+1}-a_k<M\sqrt{e^{2d}-1},\e which holds true by (6.1) if we
choose a small   $\varepsilon\in (0,1)$ such that
$$\varepsilon M(e^d-1)<M\sqrt{e^{2d}-1}.$$
Next, we want to prove (6.4). Since $a_k<a_{k+1}$, we have
$$
(x_1-a_{k})^2+\sum_{i=2}^nx_i^2
=M^2e^{2d}+(a_{k+1}-a_k)(2x_1-a_k-a_{k+1})  < M^2e^{2d}, $$ which is
the second inequality in (6.4). The first inequality in (6.4) is
reduced to   $$(x_1-a_{k})^2+\sum_{i=2}^nx_i^2
=M^2e^{2d}+(a_{k+1}-a_k)(2x_1-a_k-a_{k+1})  > M^2 ,$$ which is
equivalent to \be\la{6.7}
x_1>\f{1}{2}[a_k+a_{k+1}+\f{M^2(1-e^{2d})}{a_{k+1}-a_k}].\e By the
definition of set $\{\vec{x}\in C_N(M)|\
|\vec{x}-\vec{x}_{k+1}|=Me^d, x_1<a_{k+1}\}$, we have
$$x_1\geq a_{k+1}-Me^d.$$ Therefore, in order to prove (6.7), it is enough to
show \ba
a_{k+1}-Me^d>\f{1}{2}[a_k+a_{k+1}+\f{M^2(1-e^{2d})}{a_{k+1}-a_k}],\nonumber\ea
which is equivalent to
\ba(a_{k+1}-a_k)^2-2Me^d(a_{k+1}-a_k)-M^2(1-e^{2d})>0,\nonumber\ea
i.e., \ba a_{k+1}-a_k>M(e^d+1)\ \ \hbox{or}\ \
a_{k+1}-a_k<M(e^d-1).\nonumber\ea The  last inequality is obvious
  by (6.1).  Thus, the lemma is completed.

  $\hfill \Box$

{\bf Lemma A.2} {\sl Suppose  that $b_1>0$ , $\theta \in (0,
\frac{\pi}{2})$,   $d\in (0, \f{c(H^*,\eta)}{\mu})$ such that
$1-e^d\sin\theta>0$, and   $\delta_0$ satisfies
\be1<\delta_0<\f{1-\sin\theta}{1-e^d\sin\theta}.\e Let
 $b_k=\delta_0^{k-1}b_1$, $L_k=b_k\sin\theta$,
$\vec{x}_k= (b_k, 0, \cdots , 0)$   and
$$\tilde{A}(\vec{x}_k)=\{\vec{x}=(x_1, x_2, \cdots , x_n) \ \in\ C(\theta)\ |\
L_k< |\vec{x}-\vec{x}_k|<L_ke^d,\ x_1<b_k,
  \} $$ for $k=1, 2, \cdots .$  Then the part of the boundary of
$\tilde{A}(\vec{x}_k)$,
$$S_k:=\{\vec{x}=(x_1, x_2, \cdots , x_n)\in C(\theta)\ |\  |\vec{x}-\vec{x}_k|=L_k,\ x_1<b_k,\ L_k=b_k\sin\theta\}$$
 is completely covered by $\tilde{A}(\vec{x}_{k+1})$. Thus,
 $$C(\theta)=\bigcup_{k}\tilde{A}(\vec{x}_{k}).$$
}
{\bf Proof } Denote
 $$T_k:=\{\vec{x}\in C(\theta)\ |\  |\vec{x}-\vec{x}_k|=L_ke^d,\ x_1<b_k,\ L_k=b_k\sin\theta\}.$$
 Obviously,  The distances
 from $(0,0,\cdots,0)$ to  $S_k$, $T_{k+1}$, $S_{k+1}$ are
  $b_k(1-\sin\theta)$, $b_{k+1}(1-e^d\sin\theta)$,
 $b_{k+1}(1-\sin\theta)$, respectively. By (6.8), we have
 $$b_{k+1}(1-e^d\sin\theta)<b_k(1-\sin\theta)<b_{k+1}(1-\sin\theta).$$
 We need only to prove that $S_k\cap T_{k+1}=\emptyset$ and $S_k\cap S_{k+1}=\emptyset$.

At first, we will show that $T_{k+1}$ does not touch $S_k$  for $x_1\leq b_k$. \\
Indeed, the expressions of $S_k$ and $T_{k+1}$ are
\ba(x_1-b_k)^2+\sum_{i=2}^nx_i^2&=&b_k^2\sin^2\theta,\nonumber\\
(x_1-b_k\delta_0)^2+\sum_{i=2}^nx_i^2&=&\delta_0^2b_k^2e^{2d}\sin^2\theta,\nonumber\ea
respectively. Suppose $S_k\cap T_{k+1}\neq\emptyset$, by
calculating, we can see that the
 coordinate on $x_1-$axis of the intersection point is
  $$x_1=\f{1}{2}b_k[1+\delta_0-\f{\sin^2\theta(e^{2d}-1)}{\delta_0-1}].$$
  We claim that
  \be\la{6.9}\f{1}{2}b_k[1+\delta_0-\f{\sin^2\theta(e^{2d}-1)}{\delta_0-1}]<
  b_k(1-\sin^2\theta).\e
In order to prove the claim (6.9), we need to prove the following
inequality,
  \be\la{6.10} \f{1-\sin\theta}{1-\sin\theta
  e^d}<\f{1}{1-\tan\theta\sqrt{e^{2d}-1}}.\e
In fact, since $0<\theta<\f{\pi}{2}$ and $1<e^d<\f{1}{\sin\theta}$,
then \ba
1-\tan\theta\sqrt{e^{2d}-1}>1-\tan\theta\sqrt{(\f{1}{\sin\theta})^2-1}=0.\nonumber
\ea
 Thus,  (6.10) is equivalent to
  \ba(1-\sin\theta)(1-\tan\theta\sqrt{e^{2d}-1})<1-\sin\theta
  e^d,\nonumber
\ea
  i.e., \be\la{6.11} e^d-1<\f{1-\sin\theta}{\cos\theta}\sqrt{e^{2d}-1}.\e
  In order to prove (6.11), it is enough to prove
  \ba \cos\theta\sqrt{e^d-1}<(1-\sin\theta)\sqrt{e^d+1},\nonumber
\ea
i.e.,
 \ba
2\sin\theta(1-\sin\theta)e^d<2(1-\sin\theta),\nonumber \ea which is
obvious since $e^d<\f{1}{\sin\theta}$. Therefore, (6.10) holds.

  It follows from (6.10) and
$1<\delta_0<\f{1-\sin\theta}{1-\sin\theta
  e^d}$ that $1<\delta_0<\f{1}{1-\tan\theta\sqrt{e^{2d}-1}}$  , which implies
   $(\f{\delta_0-1}{\delta_0})^2<\tan^2\theta(e^{2d}-1) $,
  i.e.,
 \be\la{6.12} \delta_0^2-1+\sin^2\theta-\sin^2\theta
e^{2d}\delta_0^2<2\delta_0(1-\sin^2\theta)-2(1-\sin^2\theta).
  \e
Since $b_k>0$ and $\delta_0>1$,   by (6.12), we obtain  (6.9).

However, it is obvious that the coordinate on $x_1-$axis of any
point in $S_k$ is
  larger than $b_k(1-\sin^2\theta)$. Thus,
 (6.9) can imply a contradiction with  $S_k\cap T_{k+1}\neq\emptyset$, therefore $T_{k+1}$ does not intersect $S_k$  for $x_1\leq
 b_k$.

Next, we prove that $S_k\cap S_{k+1}=\emptyset$. Write the
expressions of $S_k$ and $S_{k+1}$ as follows:
 \ba(x_1-b_k)^2+\sum_{i=2}^nx_i^2&=&b_k^2\sin^2\theta,\nonumber\\
(x_1-b_k\delta_0)^2+\sum_{i=2}^nx_i^2&=&\delta_0^2b_k^2\sin^2\theta,\nonumber\ea
respectively.
 Suppose $S_k\cap S_{k+1}\neq\emptyset$. By calculating, we  see that the coordinate on $x_1-$axis of
 the intersection point
 is $$x_1=\f{\delta_0+1}{2}b_k(1-\sin^2\theta),$$
 which is larger than $b_k(1-\sin^2\theta)$ by (6.8),while, $b_k(1-\sin^2\theta)$ is the coordinate on $x_1-$axis of the
 intersection of $S_k$ and $\partial C(\theta)$,  a contradiction! Therefore, $S_k$ does
 not intersect $S_{k+1}$ and hence $\tilde{A}(\vec{x}_{k+1})$ covers $S_k$
 completely.  The lemma has been proven.
$\hfill \Box$

\vskip 0.5cm
\section* { References}
\begin{enumerate}
\itemsep -2pt

\item [1] J. Clutterbuck, O. C. Schn$\ddot{u}$rer and F. Schulze,
Stability of translating solutions to mean curvature flow, Calc.
Var. 29(2007),281-293.

\item [2] T. H. Colding and W. P. Minicozzi II, Width and mean
curvature flow, preprint, Arxiv: 0705.3827vz, 2007.

\item [3] D. Gilbarg and N. Trudinger, Second order elliptic partial
differential equations, 2nd edn (Springer, 1983).

\item [4] C. Gui, H. Jian and H. J. Ju, Properties of translating solutions   to mean curvature
 flow,  Discrete Contin.
 Dyn. Syst. 28B(2010), 441-453.

 \item[5] G. Huisken and C. Sinestrari, Mean curvature flow
singularities for mean convex surfaces, Calc. Var. Partial Differ.
Equations, 8(1999), 1-14.

\item[6]G. Huisken and C. Sinestrari, Convexity estimates for
mean curvature flow and singularities of mean convex surfaces, Acta
Math., 183(1999), 45-70.

\item[7] H.Y.  Jian,   Translating solitons of mean curvature flow of noncompact
spacelike hypersurfaces in Minkowski space, J. Differential
Equations, 220(2006), 147-162.

\item [8] H. Y. Jian, Q. H. Liu and X. Q Chen, Convexity and
symmetry of translating solitons in mean curvature flows, Chin. Ann.
Math., 26B (2005), 413-422.

\item [9] H. Y. Jian, Y. N. Liu, Long-time existence of mean
curvature flow with external force fields, Pacific J. Math., 234
(2008), 311-324.

\item [10] Z. Jin, Existence of solutions of the prescribed
mean-curvature equation on unbounded domains, Proc. Royal Soc.
Edinb., 136A, 157-179, 2006.

\item [11] Z. Jin and K. Lancaster, Theorems of Phragm$\acute{e}$n-Lindel$\ddot{o}$f type
for quasilinear elliptic equations, J. Reine Angew. Math. 514
(1999), 165-197.

\item [12] H.J. Ju, J. Lu and H.Y. Jian, Translating solutions to
mean curvarure flow with a forcing term in Minkowski space,  Commu
Pure Appl Anal.  9(2010), 963-973.

\item [13] R. Lopez, Constant mean curvature graphs on unbounded convex
domains. J. Diff. Eqns 171 (2001), 54-62.

\item [14] T. Marquardt, Remark on the anisotropic prescribed mean
curvature equation on arbitrary domain, Math. Z., 264(2010),
507-511.

\item [15] F. Schulze, Evolution of convex hypersurfaces by powers of
the mean curvature, Math Z. 251(2005), 721-733.

\item [16] F. Schulze, Nonlinear Evolution by mean curvature and
isoperimetric inequalities, J. Differential Geom., 79(2008),
197-241.

\item [17] J. Serrin, The problem of Dirichlet for quasilinear equations with
many independent variables. Phil. Trans. R. Soc. Lond. A 264 (1969),
413¨C496.

\item [18] W. M. Sheng and X. J. Wang, Singularity profile in the
mean curvature flow, Methods Appl. Anal. 16(2009), 139-155.

\item [19] W. M. Sheng and C. Wu, On asymptotic behavior for
singularities of powers of mean curvature flow,  Chin. Ann. Math.
Ser. B 30 (2009), 51--66.

\item[20] X. J.Wang, Interior gradient estimates for mean curvature
equations,Math. Z.228 (1998), 73-81.

\item [21] X.J. Wang, Convex solutions to the mean curvature flow, arixv: math. DG/0404326, preprint.

\item [22] B. White, The nature of singularities in mean curvature
flow of mean convex surfaces, J.Amer.Math.Soc.,16(2003),123-138.
 \end{enumerate}

\end{document}